\documentclass[11pt,reqno]{amsart}

\usepackage{amsmath,amssymb,amsthm, amsfonts}   
\usepackage{bbm} 
\usepackage{graphicx}   
\usepackage{pstricks}   
\usepackage[colorlinks=true, linkcolor=blue, citecolor=magenta]{hyperref} 
\usepackage{url}
\usepackage{tasks}
\usepackage{enumerate}
\usepackage{caption,subcaption} 
\usepackage{tikz,pgfplots}

\usepackage{mathtools}      
\mathtoolsset{showonlyrefs} 

\usepackage{mathrsfs}
\usepackage{multicol}   
\usepackage{fullpage}   

\usepackage[foot]{amsaddr}

\newtheorem{thm}{Theorem}[section]

\newtheorem{lem}[thm]{Lemma}

\newtheorem*{acknowledgements*}{Acknowledgements}

\theoremstyle{definition}

\numberwithin{equation}{section}

\theoremstyle{remark}
\newtheorem{rem}[thm]{Remark}

\title{Uniform periodic counterexamples to Carleson's convergence problem with polynomial symbols
}
\author{Daniel Eceizabarrena}
\address{BCAM - Basque Center for Applied Mathematics, Bilbao, Spain
		  e-mail: {\tt deceizabarrena@bcamath.org}
		}

\author{Xueying Yu}
\address{Department of Mathematics, Oregon State University, USA
		  e-mail: {\tt xueying.yu@oregonstate.edu }
    }

\date{\today}

\begin{document}

\begin{abstract}
In Carleson's convergence problem for dispersive equations 
$i\, \partial_t u + P(D)u=0$ in the periodic setting $\mathbb T^d$, we prove that
the Sobolev exponent $d/(2(d+1))$ is necessary for any non-singular polynomial symbol $P$, including the natural powers of the Laplacian $\Delta^k$.
This is in contrast with the results known in the Euclidean case, in which for symbols $P(\xi) = |\xi|^a$ with $a > 1$
the exponent $d/(2(d+1))$ is sufficient, 
but we do not know if it is necessary. 
\end{abstract}

\maketitle


\section{Introduction}

For a function $P: \mathbb R^d \to \mathbb R$ and $D = (\partial_1, \ldots, \partial_d)$, 
consider the periodic initial value problem 
\begin{equation}\label{eq:Equation}
\left\{ \begin{array}{ll}
i \, \partial_t u  + P(D) u = 0, & x \in  \mathbb T^d, \quad t \in \mathbb R, \\
u(x,0) = f(x), & x \in  \mathbb T^d,
\end{array}
\right.
\end{equation}
whose solution we denote by 
\begin{equation}\label{eq:PeriodicSolution}
    e^{itP(D)}f(x) = \sum_{n \in \mathbb Z^d} \widehat f_n\, e^{2\pi i (n \cdot x + P(n)t)}.
\end{equation}
Carleson's problem of convergence to the initial datum consists in finding $s$ for which
\begin{equation}\label{eq:Convergence}
\lim_{t \to 0} e^{itP(D)}f(x) = f(x) \qquad \text{ for almost every } x \in \mathbb T^d, \qquad \forall f \in H^s(\mathbb T^d), 
\end{equation}
and it is well-known that \eqref{eq:Convergence} is equivalent to proving the maximal estimate
\begin{equation}\label{eq:Maximal_Estimate}
\big\lVert   \sup_{t > 0} |e^{itP(D)}f| \big\rVert_{L^2(\mathbb T^d)} \lesssim \lVert f \rVert_{H^s}, \qquad \forall f \in H^s(\mathbb T^d).
\end{equation}
Thus, up to losing the endpoint information, the objective is to determine the critical exponent
\begin{equation}\label{eq:CriticalExponentMax}
s_P(\mathbb T^d) = \inf \big\{ \, s \, : \, \text{ the maximal estimate } \eqref{eq:Maximal_Estimate} \text{ holds }  \,   \big\}. 
\end{equation}

For the Schr\"odinger equation, 
corresponding to $P(n) = |n|^2$,
the best result available is
\begin{equation}
\frac{d}{2(d+1)} \leq s_{|n|^2}(\mathbb T^d) \leq  \frac{d}{d+2}, 
\end{equation} 
obtained in \cite[Prop 1 in Section 2]{MoyuaVega2008} and \cite[Prop. 3.1 and 3.2]{CompaanLucaStaffilani2021}.
The problem is also open for more general symbols. 
In this article, we focus on polynomial symbols, 
including $P_{2k}(\xi) = |\xi|^{2k}$
that corresponds to the power Laplacian $P_{2k}(D) = (-\Delta)^k$ in the equation. 
We give a lower bound for $s_{P_{2k}}(\mathbb T^d)$ in Section~\ref{sec:Results}
and we discuss the state of the art for upper bounds in Section~\ref{sec:StateOfTheArt}.

\subsection{The results}\label{sec:Results}
The main result of this note is a counterexample to the maximal estimate \eqref{eq:Maximal_Estimate} for non-singular polynomial symbols.

\begin{thm}\label{thm:Lebesgue}
Let $P \in \mathbb Z [X_1, \ldots, X_d]$ be a non-singular\footnote{If $P \in \mathbb Z[X_1, \ldots, X_d]$ has $\deg P = k$, and $P_k$ is the homogeneous part of $P$, 
we ask $\nabla P_k(x) \neq 0$ for all $x \in \mathbb C^d \setminus \{0\}$.
} polynomial of degree at least 2.
Then,
\begin{equation}
    s_P(\mathbb T^d) \geq \frac{d}{2(d+1)}. 
\end{equation}
\end{thm}

We remark that this lower bound for the critical exponent is \textbf{independent of the degree of the polynomial}, 
which is in sharp contrast with the current results available in the Euclidean case (see discussion in Section~\ref{sec:Context}). 
Most relevantly, 
the result applies to $P(\xi) = |\xi|^{2k}$ for any $k \in \mathbb N$ that correspond to the higher order dispersive equations $i \, \partial_t u + \Delta^k u = 0$.

We also adapt the proof of Theorem~\ref{thm:Lebesgue} 
to give a counterexample to the \textit{fractal} maximal estimate
\begin{equation}\label{eq:Maximal_Estimate_Fractal}
\big\lVert   \sup_{t > 0} |e^{itP(D)}f| \big\rVert_{L^2(\mathbb T^d, d\mu)} \lesssim \lVert f \rVert_{H^s}, \qquad \forall f \in H^s(\mathbb T^d), 
\quad \forall \mu \in \mathcal M_\alpha, 
\end{equation}
where $0 < \alpha < d$ is given and $\mathcal M_\alpha$ is the set of $\alpha$-Frostman measures, 
that is, Borel measures with 
\begin{equation}
    \mu(B_r) \leq C_\mu\, r^\alpha, \qquad \forall r < 1, 
\end{equation}
where $B_r$ denotes any ball of radius $r$. 
The fractal maximal estimate is directly linked to the fractal convergence problem, for which convergence in \eqref{eq:Convergence} is strengthened to almost everywhere with respect to a given Hausdorff measure $\mathcal H^\alpha$.

As in \eqref{eq:CriticalExponentMax}, 
we look for the critical Sobolev exponent $s_P(\mathbb T^d, \alpha)$ corresponding to \eqref{eq:Maximal_Estimate_Fractal}.

%
%

\begin{thm}\label{thm:Fractal}
Let $\alpha < d$. 
For the same polynomials $P$ as in Theorem~\ref{thm:Lebesgue}, 
\begin{equation}
    s_P(\mathbb T^d, \alpha) \geq \frac{d}{2(d+1)}(d+1-\alpha). 
\end{equation}
%
\end{thm}

\subsection{Context}\label{sec:Context}
The convergence problem was proposed by Carleson in \cite{Carleson1980} for the Schr\"odinger equation in $\mathbb R^d$, 
and was essentially solved in \cite{DahlbergKenig1982,Bourgain2016,DuGuthLi2017,DuZhang2019}
by showing that the critical exponent is 
\begin{equation}\label{eq:Result_Euclidean}
s_{|\xi|^2}(\mathbb R^d) = \frac{d}{2(d+1)}. 
\end{equation}
The endpoint result remains open except in dimension $d=1$.

For the fractional counterparts $|\xi|^a$,
which correspond to the fractional Laplacians $(-\Delta)^{a/2}$ in the equation, 
we know that
\begin{equation}\label{eq:UpperBound_a}
   s_{|\xi|^a}(\mathbb R^d) \leq \frac{d}{2(d+1)}, \qquad \text{ for every } a > 1,
\end{equation}
shown in \cite{Sjolin1987,ChoKo2022,LiLiXiao2021}.
The case $a \leq 1$ differs considerably; we refer to \cite{ChoKoKohLee2023} for an account of the state of the art.
The bound \eqref{eq:UpperBound_a} is sharp when $d=1$ \cite{Sjolin1987}, 
but counterexamples are scarce in $d \geq 2$. 
One of the main obstructions is the loss of a simple algebraic structure of the symbol\footnote{Compare $\xi_1^a$ in $d=1$ and $(\xi_1^2 + \xi_2^2)^{a/2}$ in $d=2$, which has a more complicated structure, especially when $a \not\in 2\mathbb N$.}.
For example, Bourgain's proof \cite{Bourgain2016} of the lower bound in \eqref{eq:Result_Euclidean} crucially exploits 
the separation of variables of the symbol $|\xi|^2 = \xi_1^2 + \ldots + \xi_d^2$, and does not naturally adapt to $|\xi|^a$.  
To avoid these obstructions, 
An, Chu and Pierce \cite{AnChuPierce2023} adapted Bourgain's construction for the separated symbol
\begin{equation}\label{eq:Symbol_AnChuPierce}
\xi_1^k + \xi_2^k + \ldots + \xi_d^k, \qquad \qquad k \in \mathbb N \setminus \{1\}, 
\end{equation}
which corresponds to the differential operator $\partial_1^k + \partial_2^k + \ldots + \partial_d^k$. 
Using Deligne's theorem to estimate Weyl sums of order $k$, they proved 
\begin{equation}\label{eq:Result_In_Rd}
s_{\xi_1^k + \ldots + \xi_d^k}(\mathbb R^d) \geq \frac14 + \frac{d-1}{4((k-1)d+1)} 
.
\end{equation}
The same results hold for polynomial symbols of degree $k$ whose homogeneous part has the form
\begin{equation}
    \xi_1^k + Q_k(\xi'), \qquad \qquad \xi' = (\xi_2, \ldots, \xi_d), \quad Q_k \text{ non-singular},
\end{equation}
as shown by the first author and Ponce-Vanegas \cite{EceizabarrenaPonceVanegas2022}.
Chu and Pierce \cite{ChuPierce2023} extended these results to a broader family of symbols, classifying them according to
the amount of cross-terms.
All these symbols can be regarded as intermediate steps towards the natural dispersive equations $i\partial_t u + \Delta^k u = 0$, whose symbols $|\xi|^{2k}$ are polynomials but do not separate variables if $k \geq 2$.

Theorem~\ref{thm:Lebesgue}, 
which applies to all non-singular polynomial symbols including all of the above as well as $\Delta^k$, 
shows that in $\mathbb T^d$ 
this problem is irrelevant, 
giving a lower bound that is both independent of $k$ 
and better than the one in $\mathbb R^d$, since 
\begin{equation}
 \frac14 + \frac{1}{4( d(k-1)+ 1 )} <  \frac{d}{2(d+1)},
 \qquad \forall k > 2. 
\end{equation} 
We briefly discuss the causes of these differences in Section~\ref{sec:Differences_Euclidean_Periodic}.

\subsection{Differences between the counterexamples in $\mathbb R^d$ and $\mathbb T^d$}\label{sec:Differences_Euclidean_Periodic}
The fact that, unlike in \eqref{eq:Result_In_Rd}, the bounds in Theorems~\ref{thm:Lebesgue} and \ref{thm:Fractal} are independent of $k$ is mainly due to the following two reasons:
\begin{enumerate}[1.]
    \item In $\mathbb T^d$, our datum is Fourier localized around integers, and the solution is eventually reduced to an exponential sum. 
    This datum has an analogue in $\mathbb R^d$ (already present in \cite{LucaRogers2017}), and the corresponding solution can be reduced to the same exponential sum via a periodization argument that requires a certain phase to be small. 
    This imposes a restriction in the parameters of the problem that depends on $k$,  
    which is not present in $\mathbb T^d$ due to the discrete nature of the solution. 

    \item One way to avoid these restrictions in $\mathbb R^d$ is to keep this datum on $\mathbb R^{d-1}$, while using the Dahlberg-Kenig wave-packet \cite{DahlbergKenig1982} on the remaining variable, as Bourgain did in \cite{Bourgain2016}. 
    This solves the problem when $k=2$, but still imposes restrictions when $k > 2$, as shown by An, Chu and Pierce \cite{AnChuPierce2023}. 
    This is due, among other reasons, to the underlying exponential sum being $(n-1)$-dimensional, instead of $n$-dimensional. 

\end{enumerate}
A few further remarks comparing Theorems~\ref{thm:Lebesgue} and \ref{thm:Fractal} to the Euclidean counterparts 
are in place.

\begin{rem}
 
\begin{enumerate}
    \item[(a)] 
    The lower bound in Theorem~\ref{thm:Lebesgue} being independent of the degree of the polynomial is in line with the positive result \eqref{eq:UpperBound_a} in the Euclidean case, which holds for $|\xi|^a$ independently of $a > 1$. 
    It therefore looks reasonable to conjecture that the lower bound is sharp. This remains to be confirmed. 
    

    \item[(b)] 
    In the fractal case, 
    the current best result
in the Euclidean space \cite{DuZhang2019,LiLiXiao2021,EceizabarrenaPonceVanegas2022_2}\footnote{The exponent $s_{\text{sawtooth}}(\alpha)$, with a complicated but explicit expression, is given in \cite{EceizabarrenaPonceVanegas2022_2}. The problem is solved for $\alpha \leq d/2$, since 
$s_{|\xi|^2}(\mathbb R^d, \alpha) = (d-\alpha)/2$
when $\alpha \leq d/2$ \cite{BarceloBennetCarberyRogers2011,Zubrinic2002}.
}
is
\begin{equation}\label{eq:Fractal_Euclidean}
    s_{\text{sawtooth}}(\alpha)
    \leq s_{|\xi|^a}(\mathbb R^d, \alpha) \leq \frac{d}{2(d+1)}(d+1-\alpha),
     \qquad \text{ when } \quad  \frac{d}{2} < \alpha < d.
\end{equation} 
In the periodic case, for the Schr\"odinger equation we have  \cite{EceizabarrenaLuca2022}
\begin{equation}\label{eq:Fractal_Results_Periodic}
   \frac{d}{2(d+1)}(d+1-\alpha) \leq  s_{|\xi|^2}(\mathbb T^d, \alpha) \leq \frac{d}{2(d+2)}(d+2-\alpha), 
   \qquad \text{ when } \quad  \alpha \leq d,
\end{equation}
so the Euclidean upper bound in \eqref{eq:Fractal_Euclidean} is the best possible in $\mathbb T^d$.
By Theorem~\ref{thm:Fractal}, the lower bound in \eqref{eq:Fractal_Results_Periodic} holds for every polynomial symbol, in particular for $|\xi|^{2k}$ with $k \in \mathbb N$.

\end{enumerate}
\end{rem}

\subsection{State of the art in $\mathbb T^d$}\label{sec:StateOfTheArt}
We conclude the introduction with a brief discussion of the state of the art for the fractional symbols $|\xi|^a$ in $\mathbb T^d$, 
which we did not find explicitly in the literature.

The current best upper bounds follow from the available Strichartz estimates (see Appendix~\ref{sec:Positive_Result_Details} for a brief explanation).
The most general periodic Strichartz estimates for symbols $|\xi|^a$ are by Schippa \cite{Schippa2020} and imply
\begin{equation}\label{eq:General_Upper_Bound}
    s_{|\xi|^a}(\mathbb T^d)
    \leq 
    \left\{
    \begin{array}{ll}
       \frac{d}{d+2},   & 1 < a \leq 2, \\
       \frac{d}{d+2}\, \frac{a}{2},   & 2 < a \leq d+2, \\
       \frac{d}{2},     & a > d+2.
    \end{array}
    \right.    
\end{equation}
Notice the large gap with the lower bound $d/(2(d+1))$ given in Theorem~\ref{thm:Lebesgue} for $a \in 2\mathbb N$,
so the problem is wide open. 
For $d=1$ and $a = k \in \mathbb N$,
due to better Strichartz estimates we have
\begin{equation}
    \frac14 
    \leq 
    s_{|\xi|^k}(\mathbb T)
    \leq 
    \left\{
    \begin{array}{lllll}
       2/5,   & k = 3,  & & & \text{ Hughes and Wooley \cite{HughesWooley2021}},  \\
       \frac{1}{2} - \frac{1}{k(k+1)}, & k \geq 4, & & & \text{ Lai and Ding \cite{LaiDing2018}}.
    \end{array}
    \right.    
\end{equation}
Regarding the fractal case, \eqref{eq:General_Upper_Bound} and the transference result\footnote{This is written for the Schr\"odinger solution $e^{it\Delta}f$, but it holds with the same proof for $e^{itD^a}f$.} in \cite[Prop. 5.2]{EceizabarrenaLuca2022} imply
\begin{equation}\label{eq:General_Upper_Bound_Fractal}
    s_{|\xi|^a}(\mathbb T^d,\alpha)
    \leq 
    \left\{
    \begin{array}{ll}
       \frac{d}{2(d+2)}\, (d+2-\alpha),   & 1 < a \leq 2, \\
       \frac{d}{2(d+2)}\, (d+a-\alpha),   & 2 < a \leq d+2, \\
       \frac{d}{2(d+2)}\, (2d+2-\alpha),     & a > d+2,
    \end{array}
    \right. 
\end{equation}
and similarly, for $d=1$, 
\begin{equation}\label{eq:General_Upper_Bound_Fractal_1D}
    s_{|\xi|^k}(\mathbb T,\alpha)
    \leq 
    \left\{
    \begin{array}{ll}
       \frac12 - \frac{\alpha}{10},   & k=3, \\
       \frac{1}{2} - \frac{\alpha}{k(k+1)},     & k \geq 4.
    \end{array}
    \right. 
\end{equation}

\section{Proof of Theorem~\ref{thm:Lebesgue}}

To prove Theorem~\ref{thm:Lebesgue} we construct a counterexample to the maximal estimate \eqref{eq:Maximal_Estimate}
inspired by the work of the first author and Lucà \cite{EceizabarrenaLuca2022}, 
and of Compaan, Lucà and Staffilani \cite{CompaanLucaStaffilani2021}, 
and we combine them with techniques by An, Chu and Pierce in \cite{AnChuPierce2023} and the first author and Ponce-Vanegas in \cite{EceizabarrenaPonceVanegas2022}.

\subsection{The initial datum}\label{sec:Datum}

Let $\psi:\mathbb R \to \mathbb R$ be a smooth function with support in $[1/4,2]$ such that $\psi(\xi) = 1$ in $\xi \in [1/2, 1]$. 
Define $\varphi:\mathbb R^d \to \mathbb R$ by 
\begin{equation}
    \varphi(\xi) = \psi(\xi_1)\, \ldots \, \psi(\xi_d), \qquad \xi = (\xi_1, \ldots, \xi_d) \in \mathbb R^d. 
\end{equation}
Let $N \gg 1$
and define the datum 
\begin{equation}\label{eq:Initial_Datum}
f_N(x) = \sum_{ n \in \mathbb Z^d } \varphi\Big( \frac{n}{N} \Big) e^{2 \pi  i n \cdot x},
\end{equation}
where $n = (n_1, \ldots, n_d) \in \mathbb Z^d$ and $x = (x_1, \ldots, x_d) \in \mathbb T^d$.
Direct computation shows that
\begin{equation}\label{eq:Norm_Of_Datum}
\lVert f_N \rVert_{H^s}^2 = \sum_{ n \in \mathbb Z^d } (1 + |n|^2)^s \, \varphi\Big( \frac{n}{N} \Big)^2
\simeq 
N^{2s + d}, 
\end{equation}
while the corresponding solution has the form
\begin{equation}\label{eq:Solution}
    e^{itP(D)}f_N(x) = \sum_{n \in \mathbb Z^d} \varphi\Big( \frac{n}{N} \Big) e^{2 \pi  i (n \cdot x + P(n) t)}. 
\end{equation}

\subsection{A pointwise lower bound for the solution}
Following Bourgain's idea in \cite{Bourgain2016}, 
we evaluate the solution \eqref{eq:Solution} around rational points to ensure some constructive interference.
Let $b \in \mathbb Z^d$ and let $q \in \mathbb N$ be prime\footnote{The estimate for the Weil sum in Theorem~\ref{thm:Estimate_Weil_Sums} is the only place where we need $q$ to be prime.} 
and such that $q \ll N$. 
Let $\delta \in \mathbb R^d$. 
Define 
\begin{equation}\label{eq:Choice_x_t}
    x = \frac{b}{q} + \delta, \qquad t = \frac{1}{q}, \qquad \text{ with } \quad|\delta| \leq \frac{c}{N},
\end{equation}
with $c >0$ small enough. 
Then, 
\begin{equation}
 e^{itP(D)}f_N(x) 
 = \sum_{n \in \mathbb Z^d} \varphi\Big( \frac{n}{N} \Big) e^{2 \pi  i \frac{b \cdot n + P(n)}{q}} \, e^{2\pi i \delta \cdot n} 
 = \sum_{n \in \mathbb Z^d} \zeta(n)\,  e^{2 \pi  i \frac{b \cdot n + P(n)}{q}},
\end{equation}
where 
\begin{equation}\label{eq:Def_Zeta}
    \zeta(n) = \varphi(n/N) e^{2\pi i \delta \cdot n}.
\end{equation}
The main idea is that, 
due to the choice of $\delta$ in \eqref{eq:Choice_x_t}, $\zeta(n) \simeq 1$ when $n_i \simeq N$ for all $i=1, \ldots, d$, and zero otherwise. 
Hence, the solution behaves essentially as the Weyl-type exponential sum 
\begin{equation}
\sum_{n_i \simeq N, \,  \forall i} e^{2 \pi  i \frac{b \cdot n + P(n)}{q}}, 
\end{equation}
which will be easier to estimate. 

To materialize this, we split the sum modulo $q$. 
Write $n = mq+r$ with $m \in \mathbb Z^d$ and $r \in \mathbb F_q^d$, where by $\mathbb F_q$ we denote the field of integers modulo $q$. 
By periodicity, we have
\begin{equation}\label{eq:Solution_in_Z}
   e^{itP(D)}f_N(x)  
   =  \sum_{r \in \mathbb F_q^d} \Big( \sum_{m \in \mathbb Z^d}   \zeta(mq+r) \Big)  e^{2 \pi  i \frac{b \cdot r + P(r)}{q}}
   = \sum_{r \in \mathbb F_q^d} Z(r)  e^{2 \pi  i \frac{b \cdot r + P(r)}{q}},
\end{equation}
where we define 
\begin{equation}\label{eq:Def_Z}
    Z(r) = \sum_{m \in \mathbb Z^d}   \zeta(mq+r), \qquad  \qquad \text{ for } r \in \mathbb F_q^d. 
\end{equation}
Decompose this function into frequencies by the discrete Fourier transform
\begin{equation}\label{eq:DiscreteFT}
    Z(r) = \sum_{\ell \in \mathbb F_q^d} \widehat Z (\ell) \, e^{2\pi i \frac{r \cdot \ell}{q}}, 
    \qquad \text{ with } \qquad 
    \widehat Z (\ell) = \frac{1}{q^d} \,  \sum_{r  \in \mathbb F_q^d} Z(r) \, e^{-2\pi i \frac{r \cdot \ell}{q}}, 
\end{equation}
so that from \eqref{eq:Solution_in_Z} we get a main component and an error term, 
\begin{equation}
\begin{split}
    e^{itP(D)}f_N(x)
    & = \widehat Z(0) \,  \sum_{r \in \mathbb F_q^d}  e^{2 \pi  i \frac{b \cdot r + P(r)}{q}} 
    + 
    \sum_{r \in \mathbb F_q^d} \sum_{\ell \in \mathbb F_q^d \setminus \{ 0 \} } \widehat Z(\ell)  \,  e^{2 \pi  i \frac{(b + \ell) \cdot r + P(r)}{q}} \\
    & = M + E, 
    \end{split}
\end{equation}
of which the dominant will be the main term $M$. 
To see that, we estimate each of them separately.

\subsection*{Main term}
From the definitions of $Z$ and $\zeta$ in \eqref{eq:Def_Z}, direct computation shows 
\begin{equation}\label{eq:Z_Hat_0}
    \widehat Z(0) 
    = \frac{1}{q^d} \,  \sum_{r  \in \mathbb F_q^d} Z(r) 
    = \frac{1}{q^d} \,  \sum_{r  \in \mathbb F_q^d} \sum_{m \in \mathbb Z^d}  \zeta(mq+r)
    = \frac{1}{q^d} \, \sum_{n \in \mathbb Z^d}  \varphi\Big( \frac{n}{N} \Big) \, e^{2\pi i \delta \cdot n} 
    \simeq \frac{N^d}{q^d}. 
\end{equation}
On the other hand, 
to estimate the exponential sum, we have the following theorem
that follows from a theorem of Deligne (see \cite[Theorem 8.4]{Deligne1974} and \cite[Theorem 11.43]{IwaniecKowalski2004})
and Hilbert's Nullstellensatz as explained in \cite[Corollary 3.2]{EceizabarrenaPonceVanegas2022}.
\begin{thm}\label{thm:Estimate_Weil_Sums}
    Let $P \in \mathbb Z[X_1, \ldots, X_d]$ such that $\deg P = k$. 
    Let $P_k$ be the homogeneous part of $P$. 
    Suppose that $\nabla P_k(x) \neq 0$ for all $x \in \mathbb C^d \setminus \{0\}$. 
    Then, for every $q$ prime such that $q \nmid k$, 
    \begin{equation}
        \Big| \sum_{n \in \mathbb F_q^d} e^{2\pi i \frac{P(n)}{q}} \Big| \leq (k-1) q^{d/2}.
    \end{equation}
\end{thm}  

Based on an argument in \cite[Section 2.1]{AnChuPierce2023}, 
we can get a lower bound for the exponential sum in $M$ from Theorem~\ref{thm:Estimate_Weil_Sums} as follows. Call
\begin{equation}
    S(b) = \sum_{r \in \mathbb F^d_q} e^{2\pi i \frac{P(r) + b \cdot r}{q}}, 
    \qquad b \in \mathbb F_q^d. 
\end{equation}
Then, 
\begin{equation}
    \sum_{b \in \mathbb F_q^d} |S(b)|^2 
    = \sum_{r,s \in \mathbb F_q^d} \Big(  e^{2\pi i \frac{P(r) - P(s)}{q}}\, \sum_{b \in \mathbb F_q^d} e^{2\pi i \frac{r-s}{q} \cdot b} \Big)
    = \sum_{r \in \mathbb F_q^d} q^d = q^{2d}, 
\end{equation}
because the inner sum is zero unless $r=s$. 
Let $c > 0$ and define $G(q) = \{ \, b \in \mathbb F_q^d \, : \, |S(b)| \geq c \, q^{d/2}\,  \}$.
By Theorem~\ref{thm:Estimate_Weil_Sums} we have $|S(b)| \leq (k-1)q^{d/2}$ for all $b$, so 
\begin{equation}
    q^{2d} = \sum_{b \in \mathbb F_q^d} |S(b)|^2 \leq |G(q)| \,   (k-1)^2 q^d 
    + q^d \, ( c q^{d/2} )^2,  
\end{equation}
which implies 
\begin{equation}
    |G(q)| \geq q^d \, \frac{1 - c^2}{(k-1)^2} \simeq q^d.
\end{equation}
This shows that 
\begin{equation}\label{eq:Weil_Lower_Bound}
    \exists G(q) \subset \mathbb F_q^d \quad  \text{ with } \quad |G(q)| \gtrsim q^d 
    \qquad \text{ such that } \quad 
    \Big| \sum_{r \in \mathbb F^d_q} e^{2\pi i \frac{P(r) + b \cdot r}{q}} \Big| \gtrsim q^{d/2} \quad \text{ for all } b \in G(q).
\end{equation}
Joining this with \eqref{eq:Z_Hat_0}, we get
\begin{equation}\label{eq:Main_Term}
    M = |\widehat Z(0)| \, 
    \Big| \sum_{r \in \mathbb F_q^d}  e^{2 \pi  i \frac{b \cdot r + P(r)}{q}} \Big| 
    \gtrsim \left(\frac{N}{q} \right)^d\, q^{d/2}, 
    \qquad \forall b \in G(q). 
\end{equation}

\subsection*{Error term}
Using the Weil bound in Theorem~\ref{thm:Estimate_Weil_Sums}, 
we first bound the error term by
\begin{equation}\label{eq:Error_FirstBound}
    |E| 
    \leq \Bigg( \sup_{\ell \in \mathbb F_q^d} \Big| \sum_{r \in \mathbb F_q^d}   \,  e^{2 \pi  i \frac{(b + \ell) \cdot r + P(r)}{q}} \Big| \Bigg)  \, 
    \sum_{\ell \in \mathbb F_q^d \setminus \{ 0 \} } \big| \widehat Z(\ell) \big|
    \, \, \lesssim \, \, 
    q^{d/2}  \sum_{\ell \in \mathbb F_q^d \setminus \{ 0 \} } \big| \widehat Z(\ell) \big|.
\end{equation}
To exploit the decay of the Fourier coefficients $\widehat Z(\ell)$, we sum by parts as follows. 
Denote the canonical basis by $\{e_j\}_{j=1}^d$
and 
define the discrete second derivatives by the differences
\begin{equation}
    \Delta_j g(n) = g(n + e_j) - 2g(n) + g(n - e_j), \qquad n = (n_1, \ldots, n_d) \in \mathbb Z^d, 
\end{equation}
for any function $g:\mathbb F_q^d \to \mathbb R$. 
The discrete Laplacian is then
\begin{equation}
    \Delta g(n) = \sum_{j=1}^d \Delta_j g(n). 
\end{equation}
With these definitions, we will use the following lemma for summation by parts whose proof is elementary and we therefore omit. 

\begin{lem}\label{thm:Lemma_Sum_by_parts}
    \begin{enumerate}
        \item[(i)] For any $g,h: \mathbb F_q^d \to \mathbb R$, 
        \begin{equation}
            \sum_{r \in \mathbb F_q^d} \Delta g(r) \,  h(r) = \sum_{r \in \mathbb F_q^d}  g(r) \, \Delta h(r)
        \end{equation}
        
        \item[(ii)] For every fixed $\ell \in \mathbb F_q^d$, 
        \begin{equation}
            \Delta_j \Big( e^{-2\pi i \frac{r \cdot \ell}{q}} \Big)
            =
            -4  \sin^2 \Big(\pi \frac{\ell_j}{q}\Big) \, e^{-2\pi i \frac{r \cdot \ell}{q}},  
        \end{equation}
        where the derivative is taken in the variable $r\in \mathbb F_q^d$. Therefore, 
        \begin{equation}
            \Delta \Big( e^{-2\pi i \frac{r \cdot \ell}{q}} \Big)
            = 
            A(\ell) \, e^{-2\pi i \frac{r \cdot \ell}{q}}, 
            \qquad \text{ with } \qquad 
            A(\ell) = -4 \sum_{j=1}^d \sin^2 \big( \pi \frac{\ell_j}{q} \big).
        \end{equation}
    \end{enumerate}
\end{lem}

For $\ell \in \mathbb F_q^d \setminus \{0\}$, 
from the definition of $\widehat Z$ in \eqref{eq:DiscreteFT}
and Lemma~\ref{thm:Lemma_Sum_by_parts} (ii) we get
\begin{equation}
\begin{split}
    \widehat Z(\ell) = \frac{1}{q^d} \, \sum_{r \in \mathbb F_q^d} Z(r) \, e^{-2\pi i \frac{r \cdot \ell}{q}}
    & = \frac{1}{q^d \, A(\ell)} \, \sum_{r \in \mathbb F_q^d} Z(r) \, \Delta\Big( e^{-2\pi i \frac{r \cdot \ell}{q}} \Big) \\
    & = \frac{1}{q^d \, A(\ell)} \, \sum_{r \in \mathbb F_q^d} \Delta Z(r) \,  e^{-2\pi i \frac{r \cdot \ell}{q}},
    \end{split}
\end{equation}
and iterating,
\begin{equation}
    \widehat Z(\ell) 
    = 
    \frac{1}{q^d \, A(\ell)^\nu} \, \sum_{r \in \mathbb F_q^d} \Delta^\nu Z(r) \,  e^{-2\pi i \frac{r \cdot \ell}{q}}, \qquad \forall \nu \in \mathbb N. 
\end{equation}
Hence, 
\begin{equation}\label{eq:Bound_Z_Hat_Iterated}
    |\widehat Z(\ell)| \leq \frac{1}{q^d \, |A(\ell)|^\nu} \, \sum_{r \in \mathbb F_q^d} |\Delta^\nu Z(r)|, 
    \qquad \forall \nu \in \mathbb N.
\end{equation}
We first estimate $\Delta^\nu Z(r)$ uniformly in $r$. 
From the definition of $\zeta$ in \eqref{eq:Def_Zeta}, we get 
\begin{equation}
    |\zeta(n+e_j) - \zeta(n)| 
    \leq \sup_x |\partial_j\zeta(x)| 
    \lesssim 1/N + \delta_j \lesssim 1/N.
\end{equation}
In the same way, 
\begin{equation}
    |\Delta_j \zeta(n)| 
    \leq \sup_x |\partial_j (\zeta(x+1) - \zeta(x))|
    \leq \sup_x |\partial_j^2 \zeta(x)| 
    \lesssim \frac{1}{N^2}, 
\end{equation}
and more generally, 
\begin{equation}
        |\Delta_{j_1} \ldots \Delta_{j_\nu} \zeta(n)| \lesssim \frac{2^{2\nu}}{N^{2\nu}}, \qquad \forall \nu \in \mathbb N, 
        \qquad \forall n \in \mathbb Z^d.
\end{equation}
Hence, from the definition of $Z$ in \eqref{eq:Def_Z}, 
\begin{equation}
\begin{split}
     \Delta^\nu Z(r) 
      & = \sum_{j_1=1}^d \cdots \sum_{j_\nu=1}^d \Delta_{j_1} \cdots \Delta_{j_\nu} Z(r) \\
     & = \sum_{j_1=1}^d \cdots \sum_{j_\nu=1}^d \Big( \sum_{m \in \mathbb Z^d}\Delta_{j_1} \cdots \Delta_{j_\nu} \zeta(mq+r) \Big).
\end{split}
\end{equation}
Since $q \ll N$, from the support condition of $\zeta$ we get
\begin{equation}
     |\Delta^\nu Z(r)| \lesssim d^\nu\, \Big(\frac{N}{q}\Big)^d \, \Big(\frac{2}{N}\Big)^{2\nu}, \qquad \forall r \in \mathbb F_q^d, 
\end{equation}
so from \eqref{eq:Bound_Z_Hat_Iterated} we get 
\begin{equation}
    |\widehat Z(\ell)| \lesssim  \Big(\frac{N}{q}\Big)^d \,   \Big(\frac{4d}{N^2 |A(\ell)|}\Big)^\nu, 
    \qquad \forall \nu \in \mathbb N.
\end{equation}
Joining this with \eqref{eq:Error_FirstBound}, we get
\begin{equation}\label{eq:Bound_Of_Error_Al}
    |E| \lesssim (4d)^\nu \,  q^{d/2} \, \Big(\frac{N}{q}\Big)^d \, \frac{1}{N^{2\nu}} \sum_{\ell \in \mathbb F_q^d \setminus \{0\}}   \frac{1}{ |A(\ell)|^\nu}.
\end{equation}
Calling $\lVert x \rVert = \operatorname{dist}(x,\mathbb Z)$, we have $|\sin(\pi x)| \simeq \lVert x \rVert$, so 
\begin{equation}
    |A(\ell)| \simeq \sum_{j=1}^d \lVert \ell_j/q \rVert^2.  
\end{equation}
This function is symmetric in each component with respect to $[q/2]$, so 
\begin{equation}
    \sum_{\ell \neq 0} \frac{1}{|A(\ell)|^\nu} 
    = \sum_{\ell \in \{1, \ldots, [q/2]\}^d} \frac{2^d}{|A(\ell)|^\nu}
    \simeq 2^d q^{2\nu} \sum_{\ell \in \{1, \ldots, [q/2]\}^d} \frac{1}{|\ell|^{2\nu} } \simeq C_d \, q^{2\nu}, 
\end{equation}
by choosing any $2\nu > d$, with $C_d$ a constant that depends in the dimension $d$ only.
Plugging this in \eqref{eq:Bound_Of_Error_Al}, we get 
\begin{equation}\label{eq:Error_Term}
    |E| \lesssim_d  q^{d/2} \, \Big(\frac{N}{q}\Big)^d \, \Big( \frac{q}{N}\Big)^{2\nu}.
\end{equation}

\subsection*{The lower bound for the solution}
If estimates \eqref{eq:Main_Term} and \eqref{eq:Error_Term} hold, 
then $q \ll N$ implies
\begin{equation}\label{eq:LowerBoundSolution}
    \big| e^{itP(D)}f_N(x) \big| 
    = |M + E| 
    \simeq \Big( \frac{N}{\sqrt{q}}\Big)^d.
\end{equation}
Thus, we need to restrict $x$ and $t$ as in \eqref{eq:Choice_x_t}, 
and $b \in G(q)$ for estimate \eqref{eq:Weil_Lower_Bound} to hold. 
Actually, we need to restrict the set $G(q)$ further for technical reasons\footnote{We only need this for the fractal result in Theorem~\ref{thm:Fractal}; the proof of Theorem~\ref{thm:Lebesgue} follows with $G(q)$ unchanged.}. 
From its definition in \eqref{eq:Weil_Lower_Bound}, we remove the elements that have some zero entry and define
\begin{equation}\label{eq:Restricted_Gq}
\widetilde{G}(q) 
=  \left\{ \, (b_1, \ldots, b_d) \in G(q) \, : \, b_i \neq 0 \quad \forall i = 1, \ldots, d 
\,   \right\}.   
\end{equation}
The cardinality of this set is comparable to that of $G(q)$. 
Indeed, $\# G(q) \geq c\, q^d$ for some fixed constant $c>0$,
and the number of $(b_1, \ldots, b_d) \in \mathbb F_q^d$ that have some zero entry is $\leq d\, q^{d-1}$, 
so 
\begin{equation}
    \#\widetilde{G}(q) \geq \#G(q) - d\, q^{d-1} 
    = q^d \, \Big( c - \frac{d}{q} \Big)
    \geq \frac{c}{2}\, q^d, 
    \qquad \qquad  \text{ for } q \gg 1. 
\end{equation}
With this, we define the set of divergence to be 
\begin{equation}\label{eq:Def_X_N}
    X_N = \bigcup_{q \in \mathcal P_Q} \, \bigcup_{b \in \widetilde{G}(q)} B\Big( \frac{b}{q}, \frac{1}{N} \Big), 
    \quad \text{ where } 
    \quad Q = N^{d/(d+1)}
    \quad \text{ and }
    \quad 
    \mathcal P_Q = \{ \, q \simeq Q \, : \, q \text{ prime }  \, \}.
\end{equation}
Thus, if $x \in X_N$, then $x \in B(b/q, 1/N)$ for some $q \in \mathcal P_Q$ and $b \in \widetilde{G}(q)$, 
so choosing $t = 1/q$ both estimates \eqref{eq:Main_Term} and \eqref{eq:Error_Term} hold, hence \eqref{eq:LowerBoundSolution} too. 
Since $Q \ll N$, we conclude that for $N \gg 1$, 
\begin{equation}\label{eq:Main_Estimate}
    \forall x \in X_N, \quad \exists t \simeq \frac{1}{Q} \quad 
    \text{ such that } 
    \quad
    \Big| e^{itP(D)}f_N(x) \Big|
    \gtrsim
    \Big(\frac{N}{\sqrt{Q}}\Big)^d.
\end{equation}

\subsection{The set of divergence}
Once we have the pointwise lower bound \eqref{eq:Main_Estimate} for the solution in $X_N$, 
to disprove the maximal estimate \eqref{eq:Maximal_Estimate} it suffices to compute the measure of $X_N$. 
From $X_N \subset [0,1]^d$, we trivially have $|X_N| \leq 1$;
the choice of $Q = N^{d/(d+1)}$ is actually the smallest choice of $Q$ that maximizes this measure.
Indeed,  
\begin{equation}\label{eq:Covering_Bound}
        |X_N| 
        \leq \sum_{q \in \mathcal P_Q} \sum_{ b \,  \in \, \widetilde G(q)} \frac{1}{N^d} 
        \simeq \frac{1}{\log Q} \, \frac{Q^{d+1}}{N^d},          
    \end{equation}
which shows that if $Q \ll N^{d/(d+1)}$, then $|X_N| \ll 1$. 
On the other hand, if $Q = N^{d/(d+1)}$, 
then $|X_N| \lesssim 1/\log Q$. 
Dirichlet approximation suggests that this upper bound is attained,
since for any $x \in [0,1]^d$, there are infinitely many $q \in \mathbb N$ and $b \in \mathbb Z^d$ such that 
    \begin{equation}\label{eq:Dirichlet_Approximation}
        \Big| x_i - \frac{b_i}{q} \Big| \leq \frac{1}{q^{(d+1)/d}}, \qquad \text{ for all } i = 1, \ldots, d. 
    \end{equation}
In other words, with $N = Q^{(d+1)/d}$, 
the balls $B(b/q, 1/N)$ should suffice to cover the unit cube $[0,1]^d$; the logarithmic loss will come from the restriction of $q$ to be prime. 
We now prove this. 

\begin{lem}
    If $Q = N^{d/(d+1)}$, then $|X_N| \gtrsim 1/\log Q$.
\end{lem}
The proof follows the steps in \cite[Section 4]{AnChuPierce2023}, a sketch of which we reproduce for completeness. 
\begin{proof}
The basic idea boils down to controling the overlap among the balls in the union of $X_N$. 
For a collection of sets $\{E_j\}_{j \in J}$, define $f = \sum_j \chi_{E_j}$. Then, 
\begin{equation}\label{eq:Measure_in_L^p_norms}
    \Big( \sum_j |E_j| \Big)^2
    = \lVert f \rVert_{L^1}^2 
    \leq |\operatorname{supp}(f)| \, \lVert f \rVert_{L^2}^2
    \leq \Big| \bigcup_j E_j \Big| \, \sum_{j,k} |E_j \cap E_k|.
\end{equation}
The number of pairs of sets $(E_j,E_k)$ that overlap is at least $|J|$,
because each set overlaps itself.
Suppose that it is bounded from above by $C |J|$ for some $C > 0$. 
Then, if all the sets have comparable measures, say $|E_j| \simeq M$, 
from \eqref{eq:Measure_in_L^p_norms} we get 
\begin{equation}\label{eq:Disjointness_Condition}
    M |J| \sum_j |E_j| 
    \lesssim \Big( \sum_j |E_j| \Big)^2 
    \leq \Big| \bigcup_j E_j \Big| \, \sum_{j,k} |E_j \cap E_k|
    \lesssim  \Big| \bigcup_j E_j \Big| \, C \, |J|\,  M,
\end{equation}
which implies  
\begin{equation}\label{eq:Disjointness_Conclusion}
    \sum_j |E_j| \lesssim  \Big| \bigcup_j E_j \Big|, 
\end{equation}
so in terms of measure, the sets $\{E_j\}$ are essentially disjoint.

We use this argument in the case of $X_N$ in \eqref{eq:Def_X_N}. 
The total number of balls in the union is 
\begin{equation}
    |J| = \sum_{q \in \mathcal P_Q} \sum_{b \in \widetilde{G}(q)} 1 \simeq \frac{Q^{d+1}}{\log Q}. 
\end{equation}
On the other hand, 
two balls with centers $b/q$ and $b'/q'$ overlap if 
\begin{equation}\label{eq:Overlap_Condition}
    \Big| \frac{b_i}{q} - \frac{b'_i}{q'} \Big| \leq \frac{2}{N}, \quad \forall i = 1, \ldots, d, 
    \qquad \Longleftrightarrow \qquad 
    | b_iq' - b_i' q| \leq \frac{2\, q\, q'}{N}, \quad \forall i = 1, \ldots, d.  
\end{equation}
We count the number of such pairs with the following lemma from elementary number theory.
We include its proof in Appendix~\ref{sec:Proof_Of_NT_Lemma} for completeness. 
\begin{lem}\label{thm:Lemma_Counting}
    Let $q,q' \in \mathbb N$. The number of integer pairs $(b,b')$ that satisfy $1 \leq b \leq q$ and $1 \leq b' \leq q'$ and that 
    \begin{equation}\label{eq A}
        0 < |bq' - b'q| \leq A
    \end{equation}
    is bounded by $2A$. 
\end{lem}

With this lemma in hand, fix $q, q' \in \mathcal P_Q$. 
If $q = q'$, given that $Q \ll N$, the condition \eqref{eq:Overlap_Condition} holds only when $b_i = b_i'$ for all $i$, 
so there are $\# \widetilde{G}(q) \simeq q^d$ such pairs $(b,b')$. 
If $q \neq q'$, by Lemma~\ref{thm:Lemma_Counting} we have at most $(4qq'/N)^d$ such pairs $(b,b')$.
In all, thus, the number of intersecting pairs is at most 
\begin{equation}
    \sum_{q \in\mathcal P_Q} q^d +  \sum_{q \neq  q' \in \mathcal P_Q} \Big( \frac{4q q'}{N} \Big)^d
    \simeq Q^d \frac{Q}{\log Q} + \Big(\frac{Q^2}{N}\Big)^d \Big( \frac{Q}{\log Q} \Big)^2
    \simeq \frac{Q^{d+1}}{\log Q} \, \Big( 1 + \frac{Q^{d+1}}{N^d \log Q} \Big).
\end{equation}
Since $Q^{d+1} = N^d$, this is bounded by $Q^{d+1}/\log Q \simeq |J|$.
Hence, \eqref{eq:Disjointness_Condition} and \eqref{eq:Disjointness_Conclusion} hold, so  
\begin{equation}\label{eq:Size_X_N}
    |X_N| \gtrsim \sum_{q \in \mathcal P_Q} \sum_{b \in \widetilde{G}(q)} \Big| B\Big( \frac{b}{q}, \frac{1}{N} \Big) \Big| 
    \simeq \frac{Q^{d+1}}{\log Q} \, \frac{1}{N^d}
    \simeq \frac{1}{\log Q}. 
\end{equation}
\end{proof}

\subsection{Disproving the maximal estimate}
From \eqref{eq:Main_Estimate}, we have the following for every $N \gg 1$:
\begin{equation}
    \forall x \in X_N, \quad \exists t \simeq N^{-d/(d+1)} \quad \text{ such that } \quad \Big| e^{itP(D) } f_N(x) \Big| \gtrsim \Big( \frac{N}{\sqrt{Q}} \Big)^d. 
\end{equation}
Joining this with \eqref{eq:Norm_Of_Datum} and \eqref{eq:Size_X_N}, we get 
\begin{equation}
    \frac{ \big\lVert \, \sup_t \big| e^{itP(D) } f_N \big| \, \big\rVert_{L^2}}{\lVert f_N \rVert_{L^2}} 
    \gtrsim |X_N|^{1/2} \, \Big( \frac{N}{\sqrt{Q}} \Big)^d \, N^{-d/2}
    \simeq \frac{1}{(\log N)^{1/2}} \, N^{\frac{d}{2(d+1)} }.
\end{equation}
Hence, if $s < \frac{d}{2(d+1)}$, we get 
\begin{equation}
    \frac{\big\lVert \, \sup_t \big| e^{itP(D) } f_N \big| \, \big\rVert_{L^2}}{\lVert f_N \rVert_{H^s}} 
    \gtrsim \frac{1}{(\log N)^{1/2}} \, N^{\frac{d}{2(d+1)} - s } \to \infty, 
    \qquad \text{ as } N \to \infty, 
\end{equation}
and the maximal estimate is disproved, 
which concludes the proof of Theorem~\ref{thm:Lebesgue}.

\section{The fractal result. Proof of Theorem~\ref{thm:Fractal}}

The pointwise lower bound in \eqref{eq:Main_Estimate}
allows us to write
\begin{equation}\label{eq:Fractal_Lower_Bound}
    \frac{\Big\lVert \sup_{0<t<1} \big| e^{itP(D)}f_N \big| \Big\rVert_{L^2(d\mu)}^2}{\lVert f_N  \rVert_{H^s}^2}
    \gtrsim \left( \frac{N}{\sqrt{Q}} \right)^{2d} \, \frac{1}{N^{2s + d}} \, \mu(X_N)
\end{equation}
for a general measure $\mu$.
We need to find a sequence of $\alpha$-dimensional measures $\mu_N$
such that the right-hand side of \eqref{eq:Fractal_Lower_Bound} diverges when $N \to \infty$. 
To do that, we need to adapt the argument in the previous sections as follows:
\begin{enumerate}[1.]
    \item  Change the choice of $Q$ to
    \begin{equation}\label{eq:Fractal_Choice_Q}
        N = Q^{(d+1)/\alpha}.
    \end{equation}
    To mark the dependence of $X_N$ on $\alpha$,
    we denote $X_N = X_N^\alpha$.
    This choice of $Q$ is natural if one repeats the argument in \eqref{eq:Covering_Bound} with the scaling of the $\mathcal H^\alpha$ Hausdorff measure instead of the Lebesgue measure, 
    or if one argues like in \eqref{eq:Dirichlet_Approximation} with the Jarn\'ik-Besicovitch theorem in mind instead of the Dirichlet approximation theorem. 

    \item Instead of restricting the $L^2$ norm to $X_N$ in \eqref{eq:Fractal_Lower_Bound}, 
    we need to work with a subset $Y_N \subset X_N$ that is a disjoint union of balls. 
   
\end{enumerate}

To identify this disjoint union, we need the following lemma.
\begin{lem}\label{thm:Disjointness}
    Let the set of centers of the balls defining $X_N^\alpha$ be 
    \begin{equation}
        \mathcal C = \Big\{ \, \frac{b}{q} \, \, : \, \,  q \in \mathcal P_Q, \quad b \in \widetilde{G}(q)    \,   \Big\}
        \subset [0,1]^d \cap \mathbb Q.
    \end{equation}
    There exists a subset $\mathcal D \subset \mathcal C$
    of cardinality $\#\mathcal D \gtrsim  Q^{d+1}/\log Q$ with the property that 
    \begin{equation}
        \Big| \frac{b}{q} - \frac{b'}{q'}  \Big| \geq \frac{3}{Q^{(d+1)/d}}, 
        \qquad \text{ for all } \quad \frac{b}{q},  \frac{b'}{q'} \in \mathcal D \quad \text{ such that } \quad \frac{b}{q} \neq  \frac{b'}{q'} . 
    \end{equation}
\end{lem}

\begin{proof}
Let $\mathcal D$ be a set of rationals $b/q$ satisfying the following four properties:
\begin{enumerate}[1.]
    
    \item \[ q \in \mathcal P_Q \quad \text{ and } \quad  b \in \widetilde{G}(q), \]

    \item 
    \begin{equation}
        \operatorname{dist}\Big( \, \frac{b}{q}, \,  \mathbb R^d \setminus [0,1]^d \, \Big) > \frac{1}{Q^{(d+1)/d}},
        \qquad \text{ for all } \frac{b}{q} \in \mathcal D, 
    \end{equation}

    \item 
    \begin{equation}
        \Big| \frac{b}{q} - \frac{b'}{q'} \Big| > \frac{3}{Q^{(d+1)/d}}, 
        \qquad \text{for all different } \quad \frac{b}{q}, \, \frac{b'}{q'} \in \mathcal D,
    \end{equation}

    \item $\mathcal D$ is maximal, in the sense that any other set $\mathcal D'$ satisfying properties 1, 2 and 3 has $\#\mathcal D' \leq \#\mathcal D$. 

\end{enumerate}

Our objective is to prove that
\begin{equation}\label{eq:Cover_Disjoint}
    X_N^d = \bigcup_{b/q \in \mathcal C} B\Big( \frac{b}{q}, \frac{1}{Q^{(d+1)/d}} \Big)  
    \subset 
    \bigcup_{b/q \in \mathcal D} B\Big( \frac{b}{q}, \frac{4}{Q^{(d+1)/d}} \Big), 
\end{equation}
since combining this with the fact that $|X_N^d| \gtrsim 1/\log N $ in \eqref{eq:Size_X_N}, 
we deduce that 
\begin{equation}
    \frac{1}{\log N} \lesssim \#\mathcal D \, \frac{4^d}{Q^{d+1}},
\end{equation}
which implies $\#\mathcal D \gtrsim Q^{d+1}/ \log N \simeq Q^{d+1}/ \log Q$ as we wanted to show.

To prove \eqref{eq:Cover_Disjoint}, take $x \in X_N^d$,
so there are $q(x) \in \mathcal P_Q$ and $b(x) \in \widetilde{G}(q)$ such that 
\[x \in B\Big( \frac{b(x)}{q(x)}, \frac{1}{Q^{(d+1)/d}} \Big).\] 
First observe that from the definition of $\widetilde{G}(q)$ we have 
\begin{equation}
    \frac{1}{q} \leq \frac{b_i(x)}{q(x)} \leq \frac{q-1}{q}, \quad \forall i = 1, \ldots, d, 
\end{equation}
and hence 
\[
\operatorname{dist} \Big(\,  \frac{b(x)}{q(x)}, \,  \mathbb R^d \setminus [0,1]^d \, \Big) 
\gtrsim \frac{1}{Q} > \frac{1}{Q^{(d+1)/d}}.
\]
We now have two alternatives: 
\begin{enumerate}
    \item If $b(x)/q(x) \in \mathcal D$, there is nothing to do. 

    \item If $b(x)/q(x) \not\in \mathcal D$, then there exists $b'/q' \in \mathcal D$ such that 
    \begin{equation}
        \Big| \frac{b(x)}{q(x)} - \frac{b'}{q'} \Big| \leq \frac{3}{Q^{(d+1)/d}}.
    \end{equation}
    Indeed, if that is not the case, then the set $\mathcal D \cup \{ b(x)/q(x)\}$ would satisfy properties 1, 2 and 3 above, but then $\mathcal D$ is not maximal.  Thus, 
    \begin{equation}
        \Big| x - \frac{b'}{q'} \Big| 
        \leq \Big| x - \frac{b(x)}{q(x)} \Big| + \Big| \frac{b(x)}{q(x)} - \frac{b'}{q'} \Big|
        \leq \frac{1}{Q^{(d+1)/d}} + \frac{3}{Q^{(d+1)/d}}
        = \frac{4}{Q^{(d+1)/d}}.
    \end{equation}
\end{enumerate}
In both cases, we get 
\[x \in B\Big( \frac{b}{q}, \frac{4}{Q^{(d+1)/d}} \Big)
\qquad \text{ for some } \quad \frac{b}{q} \in \mathcal D. 
\] 
which is what we wanted to prove.
\end{proof}

With Lemma~\ref{thm:Disjointness} in hand, we define
\begin{equation}
    Y_N^\alpha = \bigcup_{b/q \in \mathcal D} B\Big( \frac{b}{q}, \frac{1}{N}  \Big), 
    \qquad N = Q^{(d+1)/\alpha}, 
\end{equation}
which is a disjoint union independently of the choice of $\alpha \leq d$, and hence
\begin{equation}
    | Y_N^\alpha | 
    = \#\mathcal D \, \frac{1}{N^d}
    \simeq \frac{Q^{d+1}}{\log Q} \, \frac{1}{N^d}
    = \frac{1}{N^{d-\alpha} \, \log Q}
    \simeq \frac{1}{N^{d-\alpha} \, \log N}.
\end{equation}
Define the measure 
\begin{equation}\label{eq:Mu_N}
    \mu_N(A) 
    = N^{d-\alpha} \, |A \cap Y_N^\alpha|, 
    \qquad \text{ for any } A \subset [0,1]^d,  
\end{equation}
the natural choice for a measure supported on $[0,1]^d$ that captures the fractal scaling of $Y_N^\alpha$.
If we prove that $\mu_N$ is an $\alpha$-Frostman measure independently of $N$, 
then we conclude the proof of Theorem~\ref{thm:Fractal}, 
since from \eqref{eq:Fractal_Lower_Bound}, 
$Y_N^\alpha \subset X_N^\alpha$, 
and the choice of $Q$ in \eqref{eq:Fractal_Choice_Q}, 
we get
\begin{equation}
\begin{split}
    \frac{\Big\lVert \sup_{0<t<1} \big| e^{itP(D)}f_N \big| \Big\rVert_{L^2(d\mu_N)}^2}{\lVert f_N  \rVert_{H^s}^2}
    & \gtrsim \left( \frac{N}{\sqrt{Q}} \right)^{2d} \, \frac{1}{N^{2s + d}} \, \mu_N(Y_N^\alpha) \\
    & \simeq  N^{\frac{d}{d+1}(d + 1 -\alpha) - 2s} / \log N,
    \end{split}
\end{equation}
which tends to infinity as $N \to \infty$ when $s < \frac{d}{2(d+1)}(d+1-\alpha)$. 
Thus, it suffices to prove the following lemma. 
\begin{lem}
    The measure $\mu_N$ defined in \eqref{eq:Mu_N} satisfies 
    \begin{equation}
        \mu_N \in \mathcal M_\alpha, \quad \text{ and } \quad 
        \langle \mu_N \rangle_\alpha \leq 1, \qquad \forall N \in \mathbb N. 
    \end{equation}
\end{lem}
\begin{proof}
We need to prove $\mu_N(B_r) \leq r^\alpha$ for every ball $B_r$ of radius $r \in (0,1)$.
We separate cases depending on $r$.
\begin{itemize}
    \item  If $r \leq 1/N$, then $B_r$ can be completely contained in a ball of $Y_N^\alpha$, so 
    \[
    \mu_N(B_r) \leq N^{d-\alpha} \, r^d
    = (rN)^{d-\alpha} \, r^\alpha 
    \leq r^\alpha. 
    \]

    \item If $1/N \leq r \leq 1/Q^{(d+1)/d}$, then $B_r$ can intersect at most one ball of $Y_N^\alpha$, and it could totally cover it. Hence, 
    \[
    \mu_N(B_r) \leq N^{d-\alpha} \, \frac{1}{N^d}
    = \frac{1}{N^\alpha}
    \leq r^\alpha. 
    \]

    \item If $1/Q^{(d+1)/d} \leq r \leq 1$, then $B_r$ may intersect at most $Q^{d+1}\, r^d$ balls of $Y_N^\alpha$, so
    \[
    \mu_N(B_r) \leq N^{d-\alpha}\, Q^{d+1} \, r^d\, \frac{1}{N^d}
    = \frac{Q^{d+1}}{N^\alpha} r^d
    = r^d \leq r^\alpha. 
    \]
\end{itemize}
\end{proof}

\appendix

\section{On positive results for the periodic fractional Schr\"odinger equation}\label{sec:Positive_Result_Details}

We reproduce the standard procedure to obtain maximal estimates from Strichartz estimates of the form 
\begin{equation}\label{eq:GenericStrichartz}
    \lVert e^{itD^a} f_N \rVert_{L^p_{x,t}} \lesssim N^{S} \, \lVert f_N \rVert_{L^2}, 
    \qquad \text{ for } f_N \text{ with } \operatorname{supp}f_N \subset B_N
\end{equation}
using the property
\begin{equation}\label{eq:LeeLemma}
        \sup_{t\in [a,b]} |\phi(t)| \lesssim |\phi(a)| + \mu^{1/p - 1} \lVert \phi' \rVert_{L^p([a,b])} + \mu^{1/p} \lVert \phi \rVert_{L^p([a,b])}
\end{equation}
which holds for $\phi$ smooth, $1 \leq p \leq \infty$ and $\mu > 0$ (see \cite[Lemma 2.5]{Lee2006}).
Using this inequality in the time variable and minimizing $\mu = N^a$,
we get
\begin{equation}
    \begin{split}
        \sup_{0<t<1} |e^{itD^a}f_N| 
        \lesssim 
        |f_N| + N^{a/p} \, \lVert e^{itD^a}f_N \rVert_{L^p_t([0,1])}. 
    \end{split}
\end{equation}
Taking the $L^p_x$ norm and combining \eqref{eq:GenericStrichartz} with Bernstein's inequality, we get the maximal estimate
\begin{equation}\label{eq:MaximalWhen_a<2}
    \Big\lVert \sup_{0<t<1} \big| e^{itD^a}f_N \big| \Big\rVert_{L^p( \mathbb T^d)}
    \lesssim \Big( N^{\frac{d}{2} - \frac{d}{p} + \epsilon} + N^{\frac{a}{p} + S} \Big)\, \lVert f_N \rVert_{L^2(\mathbb T^d)}.
\end{equation}
The results in Section~\ref{sec:StateOfTheArt} then follow from the best available Strichartz estimates of the type \eqref{eq:GenericStrichartz}. 
Indeed, based on the Bourgain-Demeter $\ell^2$-decoupling theorem, Schippa \cite{Schippa2020} proved 
\begin{equation}\label{eq:StrichartzEstimatesFractionalBest}
    \frac{\left\lVert e^{itD^a} f_N  \right\rVert_{L^p(\mathbb T \times \mathbb T^d)}}{\lVert f_N \rVert_{L^2(\mathbb T^d)}} 
    \lesssim 
    \left\{
    \begin{array}{ll}
      N^{\frac{d}{2} - \frac{d+a}{p} + \epsilon},   &  1 < a < 2, \\
      N^{\frac{d}{2} - \frac{d+2}{p} + \epsilon}   &  a > 2, 
    \end{array}
    \right.
    \qquad \text{ for } p \geq \frac{2(d+2)}{d}.
\end{equation}    
For the particular case $d=1$ and $a = k \in \mathbb N$, from Hughes-Wooley \cite{HughesWooley2021} and Lai-Ding \cite{LaiDing2018} we have
\begin{equation}\label{eq:StrichartzEstimatesFractionalBestOneDimension}
    \frac{\left\lVert e^{itD^k} f_N  \right\rVert_{L^p(\mathbb T \times \mathbb T)}}{\lVert f_N \rVert_{L^2(\mathbb T)}} 
    \lesssim 
    \left\{
    \begin{array}{lll}
      N^{\frac{1}{2} - \frac{4}{p} + \epsilon}, & \text{ for } p \geq 10, & \text{ when } k=3, \\
      N^{\frac{d}{2} - \frac{d+2}{p} + \epsilon}, & \text{ for } p \geq k(k+1), & \text{ when } k \geq 4.
    \end{array}
    \right.
\end{equation}

\section{Proof of Lemma \ref{thm:Lemma_Counting}}\label{sec:Proof_Of_NT_Lemma}
\begin{figure}
    \centering
    \includegraphics[width=\linewidth]{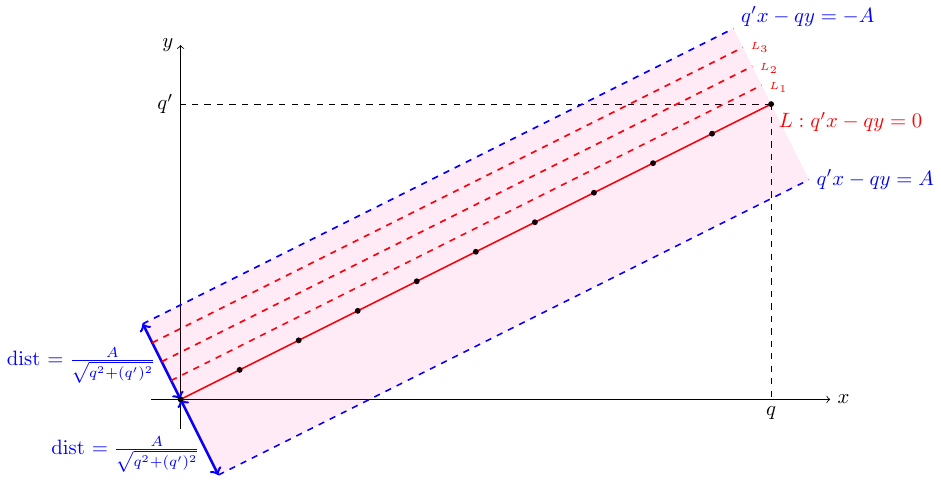}
    \caption{Lemma~\ref{thm:Lemma_Counting} counts the number of lattice points in the intersection of the diagonal rectangle with $[1,q] \times [1,q']$ that do not lie on the diagonal line $L$.}
    \label{figure}
\end{figure}

For the sake of completeness, 
we provide a geometrically motivated proof of Lemma~\ref{thm:Lemma_Counting}. 
Having $A >0$ and $q,q' \in \mathbb{N}$ fixed, and rewriting inequality \eqref{eq A} as
\begin{equation}\label{eq:Eq_Lemma_Geometric}
    0 < \frac{|bq' - b'q|}{\sqrt{q^2 + (q')^2}} \leq \frac{A}{\sqrt{q^2 + (q')^2}}, 
\end{equation}
we can interpret a pair $(b,b')$ satisfying \eqref{eq:Eq_Lemma_Geometric} as a lattice point whose distance to the line $L$ given by the equation $q' x - qy =0$ is less than $A/\sqrt{q^2 + (q')^2}$ (see Figure~\ref{figure}).
In this setting, 
it suffices to bound the number of lattice points in the diagonal rectangle shown in Figure \ref{figure} by $2A$.


Since $b,b',q,q' \in \mathbb{N}$, 
these lattice points satisfy $bq'-b'q = i$ for some $i \in \mathbb Z$ with $0 < |i| \leq A$.
Moreover, 
denoting $d = \operatorname{gcd}(q,q')$, 
we directly get $d \mid i$.
Hence, 
\begin{equation}\label{eq:Lines_Li}
    bq'-b'q = k \, d, 
    \qquad \text{ with } \quad  k \in \mathbb Z,
    \quad 
    0 < |k| \leq \lfloor A/d \rfloor,
\end{equation}
which represents a line parallel to $L$ that we denote by $L_k$.
There are at most $2A/d$ such lines. 
Hence, it suffices to bound the number of lattice points in each $L_k$ by $d$.
To do that, suppose that $(b,b') \in L_k$ is the first lattice point with $b > 0$. 
If we write $q = d \, m$ and $q' = d\, m'$ such that $\operatorname{gcd}(m,m') = 1$, 
then 
\begin{equation}
    b' = \frac{m'}{m}\, b - \frac{k}{m} \in \mathbb Z.
\end{equation}
Then, for other lattice points $(b + \beta, b' + \beta') \in \mathbb Z^2$ with $\beta \geq 0$ we have
\begin{equation}
    (b + \beta, b' + \beta') \in \mathbb Z^2 \cap L_k
    \quad \Longrightarrow \quad
    \frac{m'}{m}\, (b + \beta) - \frac{k}{m}
    \in \mathbb Z 
    \quad \Longrightarrow \quad
    \frac{m'}{m} \, \beta \in \mathbb Z,
\end{equation}
and therefore $m \mid \beta$. 
Hence, the abscissa component of the lattice points in $L_k$ have $1 \leq b + \kappa \, m \leq q$ with $0 \leq \kappa < q/m = d$. 
Consequently, each $L_k$ has at most $d$ lattice points.
Since there are at most $2A/d$ lines $L_k$, 
the maximum number of lattice points in the region is $d \times 2A/d = 2A$. 
\qed

\section*{Acknowledgements}
We would like to thank Chu-Hee Cho, Lillian Pierce, Luz Roncal and Robert Schippa for clarifying discussions.

D.E. has received funding from the European Union’s Horizon Europe research and innovation programme under Marie Sklodowska Curie Actions with grant agreement 101104250 - TIDE, 
the Simons Foundation Collaboration Grant on Wave Turbulence (Nahmod’s award ID 651469) and the American Mathematical Society and the Simons Foundation under an AMS-Simons Travel Grant for the period 2022-2024.
X.Y. is partially supported by NSF DMS-2306429.

\bibliographystyle{acm}
\bibliography{PeriodicCarlesonHomogeneous}


\end{document}